\documentclass[12pt]{article}
\usepackage{latexsym}
\usepackage{amsmath}
\usepackage{amssymb}
\usepackage{amsfonts}
\usepackage{times}
\usepackage{graphicx}
\usepackage{epsfig}
\usepackage{array}
\usepackage{flafter}
\usepackage[all]{xy}

\newtheorem{theorem}{Theorem}

\newtheorem{ExampleDef}{Example}[section]

\begin{document}
\begin{center}
{\Large {\bf Discrete Nonlinear Observers for Inertial Navigation \par}}
\vspace{1em}
{\large Yong Zhao and Jean-Jacques E. Slotine \par}
{Nonlinear Systems Laboratory \\
Massachusetts Institute of Technology \\
Cambridge, Massachusetts, 02139, USA \\
yongzhao@mit.edu, \ jjs@mit.edu \par}
\vspace{2em}
\end{center} 

\begin{abstract}
We derive an exact deterministic nonlinear observer to compute
the continuous state of an inertial navigation system based on partial
discrete measurements, the so-called strapdown problem.  Nonlinear
contraction is used as the main analysis tool, and the hierarchical
structure of the system physics is sytematically exploited. The paper 
also discusses the use of nonlinear measurements, such as distances 
to time-varying reference points.
\end{abstract}

\section{Introduction} \label{sec:introduction}

This paper derives an exact deterministic nonlinear observer to
compute the continuous state of an inertial navigation system based on
partial discrete measurements.  The main analysis tool is nonlinear
contraction theory \cite{winni98,winnithesis,winni00,chemical,jjs01}. 
Recent work on nonlinear observer design for mechanical systems based on
nonlinear contraction theory can be found in
\cite{rouchon,egeland,egeland2,jouffroy}.

Specifically, we consider the classical strap-down problem in inertial
navigation \cite{goldstein,varsha}, where angular position (Euler angles)
${\bf {x}=(\psi ,\theta ,\phi )^{T}}$ and inertial position {\bf r}
are computed from the body turn rate $\omega$ and inertial
acceleration $\gamma$, measured continuously in intrinsic (body-fixed)
coordinates,
\begin{equation} \label{eq:system}
\left\{
\begin{array}{l}
\dot{\bf x}=~{\bf H^{-1}}~{\bf \omega} \\
\dot{\bf v}=~{\bf A}~{\bf \gamma} \\
\dot{\bf r}=~~{\bf v} 
\end{array}
\right.
\end{equation} 
with
\begin{equation*}
{\bf H}=
\left [
\begin{array}{ccc}
1 & 0 & -\sin \theta \\
0 & \cos \psi & \cos \phi \sin \psi \\
0 & -\sin \psi & \cos \phi \cos \psi
\end{array} 
\right ]
\end{equation*}
and
\begin{equation*}
{\bf A}=
\left[
\begin{array}{ccc}
\cos \theta \cos \phi & \sin \psi \sin \theta \cos \phi - \cos \psi \sin \phi & \cos \psi \sin \theta \cos \phi +
 \sin \psi \sin \phi \\
\cos \theta \sin \phi & \sin \psi \sin \theta \sin \phi + \cos \psi \cos \phi & \cos \psi \sin \theta \sin \phi -
 \sin \psi \cos \phi \\
-\sin \theta & \cos \theta \sin \psi & \cos \theta \cos \psi
\end{array}
\right]
\end{equation*}
As made precise in \cite{winni00} such a system lies at
the boundary between convergence and divergence, much like a triple
integrator.  

In this paper, the continuous measurements of $\omega$ and $\gamma$
are augmented by {\it discrete} measurements of ${\bf x}$ and ${\bf
r}$, leading to a globally exponentially convergent nonlinear observer
design. Such combinations of measurements are typical in inertial
navigation, whether for vehicles or robots (see e.g. \cite{Gerdes} for
a recent discussion). The human vestibular system also features a
similar structure, with otolithic organs measuring linear acceleration
and semi-circular canals estimating angular velocity through heavily
damped angular acceleration signals, an information then combined with
visual data at much slower update rate.

Section \ref{sec:general} introduces the basic observer designs.  We
build simple observers to compute $({\bf x}, {\bf v}, {\bf r})$ based
on partial discrete measurements ${\bf x}_{i}$ and ${\bf r}_{i}$. In
Section \ref{sec:extension} we discuss extensions, such as the use of
nonlinear measurements, and the effects of system disturbance and
measurement disturbance \cite{jjs01}. We also study the case where the
inertial navigation system is expressed in quaternion form
\cite{goldstein,grassia,hestenes}. Section \ref{sec:simulation}
presents simulation results on a 3-dimensional system. Brief
concluding remarks are offered in Section \ref{sec:conclusion}.  A
very brief review of basic results in contraction theory is presented
in the Appendix.

\section{Basic Algorithm} \label{sec:general}

In this section, we construct a discrete observer for system
(\ref{eq:system}), which consists of a hierarchy of three sub-systems,
mirroring the hierarchical nature of systems physics
(\ref{eq:system}). The observer is based on the partial-measurements
of the state {\bf x} and {\bf r} at a series of instants $\{ t_{i}
\}$. \\ \\ {\bf First}, based on the discrete measurement ${\bf
x}_i$, compute {\bf x} with the observer
\begin{equation} \label{eq:observer_Ep1}
\left \{
\begin{array}{l}
\dot{\hat{\bf x}}={\bf H}^{-1}(\hat{\bf x})~{\bf \omega}  \\ \\
\hat{\bf x}^{+}_{i}=k_{1i}~\hat{\bf x}^{-}_{i}+(1-k_{1i})~{\bf x}_{i}
\end{array}
\right.
\end{equation}
where the first equation describes a continuous update between
measurements, and the second equation a discrete measurement
incorporation.

Virtual displacements, which are systematically used in mathematical
physics and optimization theory, also represent basic tools in
contraction theory (see Appendix). Computing virtual displacements in
(\ref{eq:observer_Ep1}) leads to
\begin{equation} \label{eq:delta_x_hat}
\left \{
\begin{array}{l}
\delta \dot{\hat{\bf x}}=\frac{\partial ({\bf H}^{-1}{\bf \omega})}{\partial {\hat{\bf x}}}~\delta {\hat{\bf x}} \\ \\
\delta \hat{\bf x}^{+}_{i}=k_{1i}~\delta \hat{\bf x}^{-}_{i}
\end{array}
\right.
\end{equation} 
Based on \cite{winni00}, define $\ \ \delta \hat{\bf z}={\bf \Theta}~
\delta \hat{\bf x}\ \ $ with $\ {\bf \Theta}(\hat{\bf x},t)=\bf AH\ $. 
This implies that $\ \ \ \ \left \{
\begin{array}{l}
\delta \hat{\bf z}_i^+={\bf \Theta}_i~ \delta \hat{\bf x}_i^+  \\
\delta \hat{\bf z}_i^-={\bf \Theta}_i~ \delta \hat{\bf x}_i^-
\end{array}
\right.
$ \\
From (\ref{eq:delta_x_hat}), we have
\begin{equation*} %\label{eq:delta_z_hat}
\left \{
\begin{array}{l}
\delta \dot{\hat{\bf z}}=(\dot{\bf \Theta}+{\bf \Theta}~\frac{\partial ({\bf H}^{-1}{\bf \omega})}{\partial {\hat{\bf x}}})
{\bf \Theta}^{-1}~{\delta \hat{\bf z}}={\bf 0} \\ \\
\delta \hat{\bf z}^{+}_{i}=k_{1i}~\delta \hat{\bf z}^{-}_{i}
\end{array}
\right.
\end{equation*} 
From hybrid contraction condition (\ref{eq:hybrid_condition})
in the Appendix, if
\begin{equation} \label{eq:condition_Ep1}
{\bar \lambda_{1i}}~ e^{0 \cdot \Delta t_{i}} = {\bar \lambda_{1i}} < 1 \qquad \textrm{uniformly}
\end{equation}
where ${\bar \lambda}_{1i}=k_{1i}^2$,  
then  both $\delta \hat{\bf z}$
and $\delta \hat{\bf x}$ tend to zero exponentially. So
${\hat{\bf x}}$ tends to {\bf x} exponentially.  \\ \\ {\bf Second},
 based on the discrete measurement of ${\bf r}$, compute {\bf v} with the observer
\begin{equation} \label{eq:observer_Ep2}
\left \{
\begin{array}{l}
\dot{\hat{\bf v}}={\bf A}(\hat{\bf x})~{\bf \gamma}  \\ \\
\hat{\bf v}^{+}_{i+1}=\hat{\bf v}^{-}_{i+1}-\frac {1}{\Delta t_i}~\int_{t_i}^{t_{i+1}}\hat{\bf v}dt+\frac {1}{\Delta t_i}~({\bf r}_{i+1}-{\bf r}_i)
\end{array}
\right.
\end{equation} 
From (\ref{eq:observer_Ep2}) and the first step, we get
\begin{equation} \label{eq:dp_delta_v}
\left \{
\begin{array}{l}
\frac{d}{dt}(\delta{\hat{\bf v}})=\frac{\partial{({\bf A}{\bf \gamma})}}{\partial{\hat{\bf x}}}~\delta{\hat{\bf x}}\to {\bf 0} \\ \\
\delta{\hat{\bf v}}^{+}_{i+1}=\delta {\hat{\bf v}}^{-}_{i+1}-\frac {1}{\Delta t_i} \int_{t_i}^{t_{i+1}} \delta \hat{\bf v}~dt
\end{array}
\right.
\end{equation} 
Since $\delta \hat{\bf v}$ tends exponentially to a constant, we have
$$
\frac {1}{\Delta t_i} \int_{t_i}^{t_{i+1}} \delta \hat{\bf v}~dt \to \frac {1}{\Delta t_i}~(\delta \hat{\bf v}_{i+1}^-~\Delta t_i)=\delta \hat{\bf v}_{i+1}^-
$$ 
Using (\ref{eq:dp_delta_v}), this implies that $\ \delta{\hat{\bf
v}}^{+}_{i+1}\to {\bf 0} \ $, which by continuity implies that the
constant which $\delta \hat{\bf v}$ tends to must be zero. We thus
have, exponentially,
\begin{equation} 
\left \{
\begin{array}{l}
\delta{\hat{\bf v}} \ \ \ \ \ \to \ {\bf 0} \\
\delta{\hat{\bf v}}^{+}_{i+1} \ \to \ {\bf 0}  \nonumber
\end{array}
\right.
\end{equation} 
Since by design $\hat{\bf v} = {\bf v}$ is a particular solution of
(\ref{eq:observer_Ep2}), this implies that $\hat{\bf v}$ tends to
${\bf v}$ exponentially. \\ \\ 
{\bf Third}, based on the discrete measurement ${\bf r}_i$, use the
observer
\begin{equation} \label{eq:observer_Ep3}
\left \{
\begin{array}{l}
\dot{\hat{\bf r}}=\hat{\bf v}  \\ \\
\hat{\bf r}^{+}_{i}={\bf F}_{3i}~\hat{\bf r}^{-}_{i}+({\bf I}-{\bf F}_{3i})~{\bf r}_{i}
\end{array}
\right.
\end{equation} 
Since we know $\delta \hat{\bf v}$ tends to zero exponentially, we have
$$
\left \{
\begin{array}{l}
\frac{d}{dt}(\delta{\hat{\bf r}})=\delta{\hat{\bf v}}\to {\bf 0} \\ \\
\delta{\hat{\bf r}}^{+}_{i}={\bar \lambda_{3i}}~\delta {\hat{\bf r}}^{-}_{i}
\end{array}
\right.
$$ 
If $\bar \lambda_{3i}~<~1$, i.e.
\begin{equation} \label{eq:condition_Ep3}
\bar \lambda_{3i}~e^{0 \cdot \Delta t_{i}}<1 \qquad \textrm{uniformly}
\end{equation}
where $\bar \lambda_{3i}$ is the largest eigenvalue of ${\bf F}_{3i}^{T}{\bf F}_{3i}$. So $\hat{\bf r}$ tends to {\bf r} exponentially. \\

\noindent{\bf Extension 1}: When we compute {\bf v} and {\bf r}, we only use the
discrete-time measurement ${\bf r}_i$ without ${\bf x}_i$. This allows
${\bf x}_i$ and ${\bf r}_i$ to be measured at different instants, with
the same computation. \\ \\
{\bf Extension 2}: The metric can also be written ${\bf \Theta}^T {\bf \Theta}=({\bf AH})^T({\bf AH})={\bf H}^T{\bf H}$ since {\bf A} is orthogonal.
So we can simply use ${\bf \Theta}={\bf H}$. \\ \\
{\bf Extension 3}: Assume that in (\ref{eq:observer_Ep1}) we replace the discrete update law by the more general
$$
\hat{\bf x}^{+}_{i}={\bf F}_{1i}~\hat{\bf x}^{-}_{i}+({\bf I}-{\bf F}_{1i})~{\bf x}_{i}
$$ 
where ${\bf \Theta}_i$ and ${\bf F}_{1i}$ commute. Then
\begin{equation*} 
\left \{
\begin{array}{l}
\delta \dot{\hat{\bf z}}=(\dot{\bf \Theta}+{\bf \Theta}~\frac{\partial ({\bf H}^{-1}{\bf \omega})}{\partial {\hat{\bf x}}})
{\bf \Theta}^{-1}~{\delta \hat{\bf z}}={\bf 0} \\ \\
\delta \hat{\bf z}^{+}_{i}={\bf F}_{1i}~\delta \hat{\bf z}^{-}_{i}
\end{array}
\right.
\end{equation*} 
The hybrid contraction condition (\ref{eq:condition_Ep1}) becomes
\begin{equation*}
{\bar \lambda_{1i}}~ e^{0 \cdot \Delta t_{i}} = {\bar \lambda_{1i}} < 1 \qquad \textrm{uniformly}
\end{equation*}
where ${\bar \lambda}_{1i}$ is the largest eigenvalue of ${\bf F}_{1i}^{T}{\bf F}_{1i}$.

Note that because the generalized Jacobians are zero at each step of
the hierarchy, the hybrid contraction conditions simply define the
{\it metrics} in which the discrete measurement incorporation steps
should be contracting. As we shall see later, the flexibility offered
within this constraint will allow us to trade-off model error vs
measurement error, similarly in spirit to a standard Kalman filter.

%
%section 3
%
\section{Extensions of the Basic Algorithm} \label{sec:extension}

Discussions about full discrete measurements, use of a linear
observer, nonlinear measurement, disturbance effects, and quaternion
representation are offered in this section. An observer based on full
measurement is described in Section \ref{sec:full}. Effects of system
disturbance and measurement disturbance are discussed in Section
\ref{sec:error}. Section \ref{sec:nonlinear} we develop a more general
discrete observer applicable to nonlinear measurements. Use of
quaternions is studied in Section \ref{sec:quaternion}.

%
%section 3.1
%
\subsection{Computation with Full Discrete Measurement} \label{sec:full}

Assume that {\it all} states {\bf x}, {\bf v}, and {\bf r} are actually
measured, at a series of discrete instants $\lbrace t_{i} \rbrace
$. Then steps 1 and 3 are unchanged, but we can replace step 2 (the
estimation of ${\bf v}$) by the observer
\begin{equation*} %\label{eq:observer_Ef2}
\left \{
\begin{array}{l}
\dot{\hat{\bf v}}={\bf A}(\hat{\bf x})~{\bf \gamma}  \\ \\
\hat{\bf v}^{+}_{i}={\bf F}_{2i}~\hat{\bf v}^{-}_{i}+({\bf I}-{\bf F}_{2i})~{\bf v}_{i}
\end{array}
\right.
\end{equation*} 
Since we know $\delta \hat{\bf x}$ tends to zero exponentially, we have
$$
\left \{
\begin{array}{l}
\frac{d}{dt}(\delta{\hat{\bf v}})=\frac{\partial{({\bf A}{\bf \gamma})}}{\partial{\hat{\bf x}}}~\delta{\hat{\bf x}}\to {\bf 0} \\ \\
\delta{\hat{\bf v}}^{+}_{i}={\bar \lambda_{2i}}~\delta {\hat{\bf v}}^{-}_{i}
\end{array}
\right.
$$ 
With $\bar \lambda_{2i}~<~1$, we have 
\begin{equation*} %\label{eq:condition_Ef2}
\bar \lambda_{2i}~e^{0 \cdot \Delta t_{i}}<1 \qquad \textrm{uniformly}
\end{equation*}
where $\bar \lambda_{2i}$ is the largest eigenvalue of ${\bf F}_{2i}^{T}{\bf F}_{2i}$. So $\hat{\bf v}$ tends to {\bf v} exponentially.

Note that in some cases one only needs to estimate orientation ${\bf x}$
and velocity ${\bf v}$, and that the discrete measurement of ${\bf v}$ may be
obtained from optical flow, which can be computationally "expensive"
and thus infrequent.
%
%section 3.3
%
\subsection{Disturbance Effects} \label{sec:error}

Effects of bounded inputs and measurement disturbances can be quantified
and obeserver gains chosen accordingly. 

Consider input disturbance ${\bf d}$ and measurement disturbance ${\bf
n}$, with $\|{\bf d}\| \le D$ and $\|{\bf n}\| \le N$, leading to the modified
system
\begin{equation*} %\label{eq:system_error}
\left \{
\begin{array}{l}
\dot{\hat{\bf x}}={\bf f}(\hat{\bf x})+{\bf d} \\
\hat{\bf x}_i^+=\hat{\bf x}_i^-+(k_j-1)(\hat{\bf x}_i^-+{\bf n}-{\bf x}_i)
\end{array}
\right.
\end{equation*} 

Using the basic robustness result in~\cite{winni98,jjs01}, we can quantify
the corresponding quadratic bounds $R$ on the estimation error
$$
R^{new} = \mid k_j \mid ~e^{\bar{\lambda}~\Delta t_i}~R^{old}+
\mid k_j \mid \frac{D}{\bar{\lambda}}(e^{\bar{\lambda}~\Delta t_i}-1)
 + \mid k_j-1 \mid N
$$ 
where $\bar{\lambda}$ is the largest eigenvalue of the symmetric part of $\frac{\partial {\bf f}}{\partial \hat{\bf x}}$. 

Define the objective function ($0\le k_j < 1$)
\begin{eqnarray*} %\label{eq:objective}
F(k_j)= \mid k_j \mid ~e^{\bar{\lambda}~\Delta t_i}~R^{old}+
\mid k_j \mid \frac{D}{\bar{\lambda}}(e^{\bar{\lambda}~\Delta t_i}-1)
 + \mid k_j-1 \mid N \\
= k_j ~e^{\bar{\lambda}~\Delta t_i}~R^{old}+
 k_j  \frac{D}{\bar{\lambda}}(e^{\bar{\lambda}~\Delta t_i}-1)
 + (1- k_j) N \qquad \quad 
\end{eqnarray*}
Then, $F(k_j)=(A+B-N)k_j+N$, where $A=e^{\bar{\lambda}~\Delta t_i}~R^{old}$ and
 $B=(e^{\bar{\lambda}~\Delta t_i}-1)D/{\bar{\lambda}}$.

We know $k_j$ should also satisfy
$$
k_je^{\bar{\lambda}\Delta t_i}<1\qquad \textrm{uniformly}
$$
Define $k_{max}$ as an upper bound of $k_j$. Therefore, 
$$
0\le k_j \le k_{m}
$$ where $k_{m}=min(k_{max},1)$. Finally, we obtain the minimum
of $F(k_j)$
$$
F_{min}=
\left\{
\begin{array}{l}
N,~~\textrm{when}~~k_j=0  \\
(A+B-N)k_{m}+N,~~\textrm{when}~~k_j=k_m \\
N,~~\textrm{when}~~0\le k_j \le k_m
\end{array}
\right.
\quad
\begin{array}{l}
\textrm{if}~~A+B-N>0 \\
\textrm{if}~~A+B-N<0 \\
\textrm{if}~~A+B-N=0
\end{array}
$$
where $A=e^{\bar{\lambda}~\Delta t_i}~R^{old}$ and $B=(e^{\bar{\lambda}~\Delta t_i}-1)D/{\bar{\lambda}}$.

When different measurements are available, the above formulas can also
be used to select {\it a priori} the most informative
measurement. This can be the case for instance for selecting the
direction of gaze of the eyes in hopping robot \cite{raibert}. This
can also be the case when the measurements are "expensive", for
instance computationally. \\

\noindent{\bf Extension}: The discussions above will still work when
the bounds of input disturbance and measurement disturbance are
time-varying.  If $\|{\bf d}\| \le D_i$ and $\|{\bf n}\| \le N_i$ when
$t \in \lbrack t_i,t_{i+1})$. Similar to the above, we have
$$
F_{min}=
\left\{
\begin{array}{l}
N_i,~~\textrm{when}~~k_j=0  \\
(A+B_i-N_i)k_{m}+N_i,~~\textrm{when}~~k_j=k_m \\
N_i,~~\textrm{when}~~0\le k_j \le k_m
\end{array}
\right.
\quad
\begin{array}{l}
\textrm{if}~~A+B_i-N_i>0 \\
\textrm{if}~~A+B_i-N_i<0 \\
\textrm{if}~~A+B_i-N_i=0
\end{array}
$$ 
where $A=e^{\bar{\lambda}~\Delta t_i}~R^{old}$ and $B_i=(e^{\bar{\lambda}~\Delta t_i}-1)D_i/{\bar{\lambda}}$.

%
%section 3.4
%
\subsection{Nonlinear measurements} \label{sec:nonlinear}

For the system $\dot{\bf x}={\bf f}(\bf x)$, consider the observer
\begin{equation} \label{eq:nonlinear_observer}
\left \{
\begin{array}{l}
\dot{\hat{\bf x}}={\bf f}(\hat{\bf x}) \\
\hat{\bf x}_i^+=\hat{\bf x}_i^-+{\bf g}_i(\hat{\bf y}_i^-)-{\bf g}_i({\bf y}_i)
\end{array}
\right.
\end{equation}
where
\begin{equation*}
\left \{
\begin{array}{l}
{\bf y}_i={\bf y}_i({\bf x}_i) \\
\hat{\bf y}_i^-={\bf y}_i(\hat{\bf x}_i^-)
\end{array}
\right.
\end{equation*} 
We have
\begin{equation} \label{eq:nonlinear_delta}
\left \{
\begin{array}{l}
\delta \dot{\hat{\bf x}}=\frac{\partial {\bf f}}{\partial \hat{\bf x}}~\delta \hat{\bf x} \\ \\
\delta \hat{\bf x}_i^+=({\bf I}+\frac{\partial {\bf g}_i}{\partial \hat{\bf y}_i} \frac{\partial \hat{\bf y}_i}{\partial \hat{\bf x}_i})\delta \hat{\bf x}_i^-
\end{array}
\right.
\end{equation}

Defining $\ \delta \hat{\bf z}={\bf \Theta}~ \delta \hat{\bf x}\ $, we have 
$
\left \{
\begin{array}{l}
\delta \hat{\bf z}_i^+={\bf \Theta}_i~ \delta \hat{\bf x}_i^+  \\
\delta \hat{\bf z}_i^-={\bf \Theta}_i~ \delta \hat{\bf x}_i^-
\end{array}
\right.
$. \  Using Equation (\ref{eq:nonlinear_delta}) yields
\begin{equation*} %\label{eq:delta_z_hat}
\left \{
\begin{array}{l}
\delta \dot{\hat{\bf z}}={\bf F}~{\delta \hat{\bf z}} \\ 
\delta \hat{\bf z}^{+}_{i}={\bf F}_{i}~\delta \hat{\bf z}^{-}_{i}
\end{array}
\right.
\end{equation*} 
where $
{\bf F}=(\dot{\bf \Theta}+{\bf \Theta}~\frac{\partial {\bf f}}
{\partial {\hat{\bf x}}}){\bf \Theta}^{-1}
$ and $
{\bf F}_i={\bf \Theta}_i({\bf I}+\frac{\partial {\bf g}_i}{\partial \hat{\bf y}_i} \frac{\partial \hat{\bf y}_i}{\partial \hat{\bf x}_i}){\bf \Theta}_i^{-1}
$.
The sufficient contraction condition on hybrid systems
can be written
\begin{equation} \label{eq:condition_gn}
\bar{\lambda}_i~e^{\bar{\lambda}~\Delta t_i}<1
\end{equation}
where $\bar{\lambda}_i=\lambda_{max}({\bf F}_i^T{\bf F}_i)$ and
$\bar{\lambda}$ is the largest eigenvalue of the symmetric matrix
${\bf F}^T+{\bf F}$.  If condition (\ref{eq:condition_gn}) is
satisfied by an appropriate choice of ${\bf g}_i$, then $\hat{\bf x}$
will tend to {\bf x} exponentially.

A a simple illustration, consider using distance measurements instead
of direct cartesian position measurements. In the 3-dimensional space,
measure the distances from one point $X=(x_1,x_2,x_3)^T$ to four
time-varying reference points $A=[a_1(t),a_2(t),a_3(t)]^T$,
$B=[b_1(t),b_2(t),b_3(t)]^T$, $C=[c_1(t),c_2(t),c_3(t)]^T$, and
$D=[d_1(t),d_2(t),d_3(t)]^T$,
%$A=(a_1,a_2,a_3)^T$, $B=(b_1,b_2,b_3)^T$, $C=(c_1,c_2,c_3)^T$, and $D=(d_1,d_2,d_3)^T$. 
\begin{equation} \label{eq:nonlinear_y}
\begin{array}{l}
y_1=\mid XA \mid =\sqrt{(x_1-a_1)^2+(x_2-a_2)^2+(x_3-a_3)^2}  \\
y_2=\mid XB \mid =\sqrt{(x_1-b_1)^2+(x_2-b_2)^2+(x_3-b_3)^2}  \\ 
y_3=\mid XC \mid =\sqrt{(x_1-c_1)^2+(x_2-c_2)^2+(x_3-c_3)^2}  \\
y_4=\mid XD \mid =\sqrt{(x_1-d_1)^2+(x_2-d_2)^2+(x_3-d_3)^2}  
\end{array}
\end{equation} 
The discrete-update part of observer (\ref{eq:nonlinear_observer}) can be built up as below,
\begin{equation}  \label{eq:nonlinear_od1}
\left [
\begin{array}{l}
\hat{x}_{1,i}^+ \\
\hat{x}_{2,i}^+ \\
\hat{x}_{3,i}^+
\end{array}
\right]
=
\left [
\begin{array}{l}
\hat{x}_{1,i}^- \\
\hat{x}_{2,i}^- \\
\hat{x}_{3,i}^-
\end{array}
\right]
-\frac{1}{2}
{\bf K}_i
\left[
\begin{array}{l}
(\hat{y}_{1,i}^-)^2-(\hat{y}_{2,i}^-)^2-(y_{1,i}^2-y_{2,i}^2) \\
(\hat{y}_{2,i}^-)^2-(\hat{y}_{3,i}^-)^2-(y_{2,i}^2-y_{3,i}^2) \\
(\hat{y}_{3,i}^-)^2-(\hat{y}_{4,i}^-)^2-(y_{3,i}^2-y_{4,i}^2)
\end{array}
\right]
\end{equation}
where ${\bf K}_i$ is a 3 by 3 time-varying gain matrix.
Using equation (\ref{eq:nonlinear_y}) yields
\begin{equation} \label{eq:nonlinear_dvarying}
\delta \hat{\bf x}_i^+ 
=
({\bf I}-
{\bf K}_i~{\bf J}_i)
\delta \hat{\bf x}_i^-
\end{equation}
$$
\textrm{where}
~~
{\bf J}_i=
\left[
\begin{array}{ccc}
(b_{1i}-a_{1i}) & (b_{2i}-a_{2i}) & (b_{3i}-a_{3i}) \\ 
(c_{1i}-b_{1i}) & (c_{2i}-b_{2i}) & (c_{3i}-b_{3i}) \\ 
(d_{1i}-c_{1i}) & (d_{2i}-c_{2i}) & (d_{3i}-c_{3i})
\end{array}
\right]
$$ 
where subscript $i$ refers to the value at time $t_i$. \\

Assume ${\bf J}_i$ is non-singular. Then we can choose 
\begin{equation} \label{eq:nonlinear_K_i}
{\bf K}_i=k_i~{\bf J}_i^{-1}
\end{equation}
With Equation (\ref{eq:nonlinear_dvarying}), we have 
$$
\delta \hat{\bf x}_i^+ = (1-k_i)\delta \hat{\bf x}_i^-
$$
By choosing $k_i$, we can make $\bar{\lambda}_i$ satisfy the following contraction contidion that makes $\delta \hat{\bf z}$
 tends to zero 
exponentially.\begin{equation} \label{eq:condition_timevarying}
\bar{\lambda}_i~e^{\bar{\lambda}~\Delta t_i}<1
\end{equation}
where $\bar{\lambda}_i=(1-k_i)^2$ and $\bar{\lambda}$ is the largest
eigenvalue of the symmetric matrix ${\bf F}^T+{\bf F}$.  Therefore,
$\delta \hat{\bf x}$ will tend to zero, and $\hat{\bf x}$ will tend to
{\bf x} exponentially. \\

\noindent{\bf Remark}\ \  When ${\bf J}_i$ is singular,  one has
\begin{equation} \label{eq:nonlinear_singular1}
\left|
\begin{array}{ccc}
(b_{1i}-a_{1i}) & (b_{2i}-a_{2i}) & (b_{3i}-a_{3i}) \\ 
(c_{1i}-b_{1i}) & (c_{2i}-b_{2i}) & (c_{3i}-b_{3i}) \\ 
(d_{1i}-c_{1i}) & (d_{2i}-c_{2i}) & (d_{3i}-c_{3i})
\end{array}
\right|=0
\end{equation} \\
Equation (\ref{eq:nonlinear_singular1}) is equivalent to
\begin{equation} \label{eq:nonlinear_singular2}
[(b_{1i}-a_{1i}){\bf i}+(b_{2i}-a_{2i}){\bf j}+(b_{3i}-a_{3i}){\bf k}] \cdot 
\left|
\begin{array}{ccc}
{\bf i} & {\bf j} & {\bf k} \\ 
(c_{1i}-b_{1i}) & (c_{2i}-b_{2i}) & (c_{3i}-b_{3i}) \\ 
(d_{1i}-c_{1i}) & (d_{2i}-c_{2i}) & (d_{3i}-c_{3i})
\end{array}
\right|=0
\end{equation} \\
which we can write
$$
\overrightarrow{AB} \cdot (\overrightarrow{BC} \times \overrightarrow{CD})=0
$$
This means that points A, B, C, and D are in the same plane, and therefore 
that the geometry does not contain enough information to infer position.  \\ %\ref{fig:nonlinear_plant}

To compute velocity, one can rewrite observer (\ref{eq:observer_Ep2}) as
{\scriptsize
\begin{equation} \label{eq:observer_Ep2_distance}
\left\{
\begin{array}{l}
\dot{\hat{\bf v}}={\bf A}(\hat{\bf x})~{\bf \gamma}  \\ 
\hat{\bf v}_{i+1}^+=\hat{\bf v}_{i+1}^--\frac{1}{2}\frac{1}{\Delta t_i}
\{
{\bf K}_{i+1}
\left[
\begin{array}{c}
(\hat{y}_{1,i+1}^-)^2-(\hat{y}_{2,i+1}^-)^2-(y_{1,i+1}^2-y_{2,i+1}^2) \\
(\hat{y}_{2,i+1}^-)^2-(\hat{y}_{3,i+1}^-)^2-(y_{2,i+1}^2-y_{3,i+1}^2) \\
(\hat{y}_{3,i+1}^-)^2-(\hat{y}_{4,i+1}^-)^2-(y_{3,i+1}^2-y_{4,i+1}^2)
\end{array}
\right]
-
{\bf K}_i
\left[
\begin{array}{c}
(\hat{y}_{1,i}^+)^2-(\hat{y}_{2,i}^+)^2-(y_{1,i}^2-y_{2,i}^2) \\
(\hat{y}_{2,i}^+)^2-(\hat{y}_{3,i}^+)^2-(y_{2,i}^2-y_{3,i}^2) \\
(\hat{y}_{3,i}^+)^2-(\hat{y}_{4,i}^+)^2-(y_{3,i}^2-y_{4,i}^2)
\end{array}
\right]
\}
\end{array}
\right.
\end{equation}
}
where  
\begin{equation*} %\label{eq:nonlinear_y}
\begin{array}{l}
y_1=\mid {\bf r}A \mid =\sqrt{(r_1-a_1)^2+(r_2-a_2)^2+(r_3-a_3)^2}  \\
y_2=\mid {\bf r}B \mid =\sqrt{(r_1-b_1)^2+(r_2-b_2)^2+(r_3-b_3)^2}  \\ 
y_3=\mid {\bf r}C \mid =\sqrt{(r_1-c_1)^2+(r_2-c_2)^2+(r_3-c_3)^2}  \\
y_4=\mid {\bf r}D \mid =\sqrt{(r_1-d_1)^2+(r_2-d_2)^2+(r_3-d_3))^2}  
\end{array}
\end{equation*}
and
{\scriptsize
\begin{equation} \label{eq:nonlinear_gain}
{\bf K}_i=
\left[
\begin{array}{ccc}
(b_{1i}-a_{1i}) & (b_{2i}-a_{2i}) & (b_{3i}-a_{3i}) \\ 
(c_{1i}-b_{1i}) & (c_{2i}-b_{2i}) & (c_{3i}-b_{3i}) \\ 
(d_{1i}-c_{1i}) & (d_{2i}-c_{2i}) & (d_{3i}-c_{3i})
\end{array}
\right]^{-1}
~~\textrm{and}~~
{\bf K}_{i+1}=
\left[
\begin{array}{ccc}
(b_{1i+1}-a_{1i+1}) & (b_{2i+1}-a_{2i+1}) & (b_{3i+1}-a_{3i+1}) \\ 
(c_{1i+1}-b_{1i+1}) & (c_{2i+1}-b_{2i+1}) & (c_{3i+1}-b_{3i+1}) \\ 
(d_{1i+1}-c_{1i+1}) & (d_{2i+1}-c_{2i+1}) & (d_{3i+1}-c_{3i+1})
\end{array}
\right]^{-1}
\end{equation}}
We then have
\begin{equation*} %\label{eq:dp_delta_v}
\left \{
\begin{array}{l}
\frac{d}{dt}(\delta{\hat{\bf v}})=\frac{\partial{({\bf A}{\bf \gamma})}}{\partial{\hat{\bf x}}}~\delta{\hat{\bf x}}\to {\bf 0} \\ \\
\delta{\hat{\bf v}}^{+}_{i+1}=\delta{\hat{\bf v}}^{-}_{i+1}-\frac {1}{\Delta t_i}(\delta{\hat{\bf r}}^{-}_{i+1}-
\delta{\hat{\bf r}}^{+}_{i})
=\delta {\hat{\bf v}}^{-}_{i+1}-\frac {1}{\Delta t_i} \int_{t_i}^{t_{i+1}} \delta \hat{\bf v}~dt
\end{array}
\right.
\end{equation*}
which is the same as equation (\ref{eq:dp_delta_v}). Similarly to the
second step of Section \ref{sec:general}, this shows that $\hat{\bf
v}$ tends to {\bf v} exponentially. \\

Note that the geometry problem of going from distances to positions is
solved by a dynamic system, the observer, rather than explicitly at
each instant. In general, one may also use linear measurements at some
instants and nonlinear ones at others.

Note that if a measurement is delayed, the algorithms work similarly
but the actual information is available after the delay (i.e. the
measurement is incorporated at some past time and the forward
simulation runs instantly to the current time).

Consider now, extending section \ref{sec:error}, the effect of model and measurement errors. For the system $\dot{\bf x}={\bf f}(\bf x)+{\bf d}$, with the following nonlinear observer,
\begin{equation*}
\left \{
\begin{array}{l}
\dot{\hat{\bf x}}={\bf f}(\hat{\bf x})+{\bf d} \\
\hat{\bf x}_i^+=\hat{\bf x}_i^-+{\bf g}_j(\hat{\bf y}_i^-)-{\bf g}_j({\bf y}_i-
{\bf n})
\end{array}
\right.
\end{equation*}
where$
\left \{
\begin{array}{l}
{\bf y}_i={\bf y}_i({\bf x}_i) \\
\hat{\bf y}_i^-={\bf y}_i(\hat{\bf x}_i^-)
\end{array}
\right.
$ and $
\left \{
\begin{array}{l}
{\bf d}-\textrm{model error}\quad \|{\bf d}\|<D \\
{\bf e}-\textrm{measurement error}\quad \|{\bf n}\|<N
\end{array}
\right.
$ \\ \\
So we have
\begin{equation*}
\left \{
\begin{array}{l}
\delta \dot{\hat{\bf x}}=\frac{\partial {\bf f}}{\partial \hat{\bf x}}(\delta \hat{\bf x})+\delta {\bf d} \\ \\
\delta \hat{\bf x}_i^+=({\bf I}+\frac{\partial {\bf g}_j}{\partial \hat{\bf y}_i} \frac{\partial \hat{\bf y}_i}{\partial \hat{\bf x}_i})\delta \hat{\bf x}_i^-+\frac{\partial {\bf g}_i}{\partial ({\bf y}_i-{\bf n})}\delta {\bf n}
\end{array}
\right.
\end{equation*}
We know the quadratic bounds $R$ on the estimation error
$$
R^{new}=\sqrt{\bar{\lambda}_i}~e^{\lambda_i\Delta t_i}~R^{old}+
\sqrt{\bar{\lambda}_i}~\frac{D}{\bar{\lambda}}(e^{\bar{\lambda}\Delta t_i}-1)
+\sqrt{\bar{\lambda}_{ei}}~N
$$
where $\bar{\lambda}_i=\lambda_{max}(({\bf I}+\frac{\partial {\bf g}_j}{\partial \hat{\bf y}_i} \frac{\partial \hat{\bf y}_i}{\partial \hat{\bf x}_i})^T({\bf I}+\frac{\partial {\bf g}_j}{\partial \hat{\bf y}_i} \frac{\partial \hat{\bf y}_i}{\partial \hat{\bf x}_i}))$, $\bar{\lambda}_{ei}=\lambda_{max}((\frac{\partial {\bf g}_i}{\partial ({\bf y}_i-{\bf n})})^T(\frac{\partial {\bf g}_i}{\partial ({\bf y}_i-{\bf n})}))$, and $\lambda_i$ is the largest eigenvalue of the symmetric part of $\frac{\partial {\bf f}}{\partial \hat{\bf x}}$. \\ \\
We can choose the most relevant discrete update function ${\bf g}_j$ which will best contribute to improving the estimate $\hat{\bf x}$ (i.e., to minimize $R^{new}$).

%
%section 3.5
%

\subsection{Quaternion Representation} \label{sec:quaternion}

Angular position can be expressed in quaternion form, avoiding
representation singularities \cite{goldstein,grassia}. Quaternions
express a rotation of angle $\theta$ about the unit vector ${\bf n}$
as ${\bf q}=(\cos (\theta /2), {\bf n} \sin(\theta /2))^T$.  With $\
{\bf q}=(q_0,q_1,q_2,q_3)^{T}$ the quaternion vector, this leads to
\begin{equation*} %\label{eq:quaternion}
\left\{
\begin{array}{l}
\dot{\bf q}=~\frac{1}{2}~{\bf \Omega}~{\bf q} \\
\dot{\bf v}=~{\bf A}~{\bf \gamma} \\
\dot{\bf r}=~~{\bf v} 
\end{array}
\right.
\end{equation*}
where
\begin{equation*}
{\bf \Omega}=
\left[
\begin{array}{cccc}
0 & -\omega_1 & -\omega_2 & -\omega_3 \\
\omega_1 & 0 & -\omega_3 & \omega_2 \\
\omega_2 & \omega_3 & 0 & -\omega_1 \\
\omega_3 & -\omega_2 & \omega_1 & 0 
\end{array} 
\right]
\end{equation*}
and
\begin{equation*}
{\bf A(q)}=
\left[
\begin{array}{ccc}
q_0^2+q_1^2-q_2^2-q_3^2 & 2(q_1q_2-q_0q_3) & 2(q_1q_3+q_0q_2) \\
2(q_1q_2+q_0q_3) & q_0^2-q_1^2+q_2^2-q_3^2 & 2(q_2q_3-q_0q_1) \\
2(q_1q_3-q_0q_2) & 2(q_2q_3+q_0q_1) & q_0^2-q_1^2-q_2^2+q_3^2
\end{array}
\right]
\end{equation*} 
In this representation, the fact that the dynamics of {\bf q} is
indifferent is obvious, since $\bf \Omega$ is skew-symmetric. 

The observers can be derived as earlier, simply by replacing
(\ref{eq:observer_Ep1}) by
\begin{equation*} %\label{eq:observer_Qf1}
\left \{
\begin{array}{l}
\dot{\hat{\bf q}}=\frac{1}{2}~{\bf \Omega}~{\bf q}  \\ \\
\hat{\bf q}^{+}_{i}={\bf F}_{1i}~\hat{\bf q}^{-}_{i}+({\bf I}-{\bf F}_{1i})~{\bf q}_{i}
\end{array}
\right.
\end{equation*}
based on the discrete measurements ${\bf q}_i$.  Computing virtual
displacements
\begin{equation*} %\label{eq:delta_q_hat}
\left \{
\begin{array}{l}
\delta \dot{\hat{\bf q}}={\bf \Omega}~\delta {\hat{\bf q}} \\ \\
\delta \hat{\bf q}^{+}_{i}={\bf F}_{1i}~\delta \hat{\bf q}^{-}_{i}
\end{array}
\right.
\end{equation*}
and because the dynamics of {\bf q} is indifferent, we only need
\begin{equation} \label{eq:condition_Qp1}
{\bar \lambda_{1i}}~ e^{0 \cdot \Delta t_{i}} = {\bar \lambda_{1i}} < 1 \qquad \textrm{uniformly}
\end{equation}
where ${\bar \lambda}_{1i}$ is the largest eigenvalue of ${\bf
F}_{1i}^{T}{\bf F}_{1i}$. Under Condition
(\ref{eq:condition_Qp1}), $\delta \hat{\bf q}$ tends to zero exponentially, and
${\hat{\bf q}}$ tends to {\bf q} exponentially. 

The other two steps are unchanged, with ${\bf A}(\hat{\bf x})$ being
replaced by ${\bf A}(\hat{\bf q})$. 

All the above variations and extensions can of course be combined. 

%
%section 4
%
\section{Simulation} \label{sec:simulation}
%\Example{}{

In this section, we will do a 3-dimentional simulation about system (\ref{eq:system}) based on the discrete measurement ${\bf x}_i$ and the nonlinear distance measurements $y_{1,i}$, $y_{2,i}$, $y_{3,i}$, and $y_{4,i}$, as in Section \ref{sec:nonlinear}.

Consider System (\ref{eq:system}) in the 3-dimensional case. Where
$$
{\bf \omega}=
~
\left[
\begin{array}{c}
\frac {2+\sin t}{3} \\ 
\frac {3+\cos t}{5} \\
\frac {2+\sin 2t}{3}
\end{array}
\right]~~\textrm{and}~~
{\bf \gamma}=~
\left[
\begin{array}{c}
\cos(2t) \\
\sin t \\
\frac{1+2\sin t}{3}
\end{array}
\right]
$$
Four time-varying reference points are chosen as below (all move on circular trajectories),
$$
A(a_1,a_2,a_3)~~
\left\{
\begin{array}{c}
a_1^2+a_2^2=1 \\
a_3=0
\end{array}
\right.
\qquad
B(b_1,b_2,b_3)~~
\left\{
\begin{array}{c}
(b_1-60)^2+b_2^2=1 \\
b_3=0
\end{array}
\right.
$$ 
$$
C(c_1,c_2,c_3)~~
\left\{
\begin{array}{c}
(c_2-60)^2+c_3^2=1 \\
c_1=60 
\end{array}
\right.
\qquad
D(d_1,d_2,d_3)~~
\left\{
\begin{array}{c}
(d_1-60)^2+(d_3-60)^2=1 \\
d_2=60 
\end{array}
\right.
$$

\begin{figure}[h]
\begin{center}
\epsfig{figure=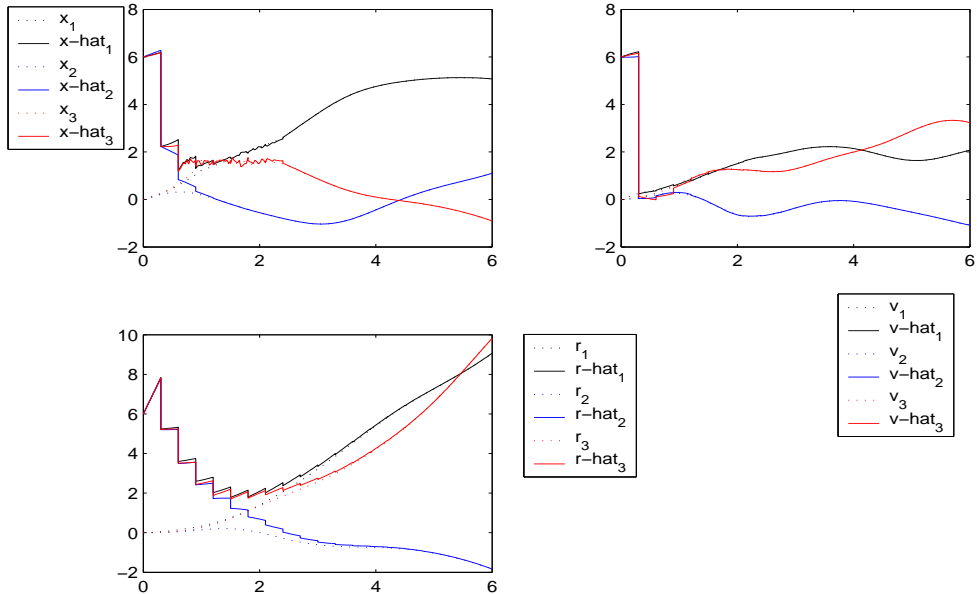,height=80mm,width=130mm}
\end{center}
\caption{\small Simulation result of computing ${\bf x}$, ${\bf v}$, and ${\bf r}$ with the discrete measurements 
${\bf x}_{i}$ and ${\bf r}_{i}$}
\label{fig:example_Ep}
\end{figure}

Observer (\ref{eq:observer_Ep1}) with $\bar \lambda_{1i}=1/9$ is used
to compute {\bf x}. Observer (\ref{eq:observer_Ep2_distance}) with
gain (\ref{eq:nonlinear_gain}) is used to compute {\bf v}. Using
observer (\ref{eq:nonlinear_observer},\ref{eq:nonlinear_od1}) and
gain (\ref{eq:nonlinear_K_i}), we choose $k_i=\frac{2}{3}$ to
satisfy Condition (\ref{eq:condition_gn}), thus we can compute {\bf
r}. Figure {\ref{fig:example_Ep}} shows $(\hat{\bf x},\hat{\bf
v},\hat{\bf r})^{T}$ tends to $({\bf x},{\bf v},{\bf r})^{T}$
exponentially.

%
%section 5
%

\section{Concluding Remarks} \label{sec:conclusion}

Observers similar to those developed in this paper can in principle be
applied to other continuous nonlinear systems besides inertial
navigation systems, although much simplification was afforded by
exploiting the hierarchical structure of the system physics. An
animation of the basic observer as applied to head stabilization
\cite{berthoz} in a simulated robot hopper \cite{raibert} can also be
found in {\em http://web.mit.edu/nsl/www/hopping\_robot.mpg}.

\noindent{\bf Acknowledgement}\ \   This paper benefited from stimulating
 discussions with Dr. Agostino Martinelli.

\renewcommand{\baselinestretch}{1.}

\section*{Appendix}
The main theorem of contraction analysis~\cite{winni98} can be stated
as

\begin{theorem} Consider the deterministic system $ \ \dot{\bf x} =
{\bf f}({\bf x},t) \ $, where ${\bf f}$ is a smooth nonlinear
function. If theres exist a uniformly positive definite metric

${\bf M(\bf x}, t) \ = \ {\bf \Theta}({\bf x}, t)^T \ {\bf
\Theta}({\bf x}, t)$

\noindent such that the associated generalized Jacobian
 
${\bf F} \ = \ \left(\dot{\bf \Theta} + {\bf \Theta} \frac{\partial {\bf
  f}} {\partial {\bf \bf x}} \right){\bf \Theta}^{-1}$ 

\noindent is uniformly negative definite, then all system trajectories
then converge exponentially to a single trajectory, with convergence
rate $| \lambda_{max} |$, where $\lambda_{max}$ is the largest
eigenvalue of the symmetric part of ${\bf F}$. The system is said to 
be contracting.
\label{th:theoremF}
\end{theorem}

It can be shown conversely that the existence of a uniformly positive
definite metric with respect to which the system is contracting is
also a necessary condition for global exponential convergence of
trajectories.  In the linear time-invariant case, a system is globally
contracting if and only if it is strictly stable, with ${\bf F}$
simply being a normal Jordan form of the system and ${\bf \Theta}$ the
coordinate transformation to that form. Furthermore, since 
\begin{small}
$$
{\bf \Theta}^{-1} \ {\bf F}_s \ {\bf \Theta} \ = \ \frac{1}{2}\ \ {\bf M}^{-1}\
(\dot{\bf M} + \ {\bf M} \ \frac{\partial {\bf f}}{\partial {\bf x}}\
+ \ \frac{\partial {\bf f}}{\partial {\bf x}}^T {\bf M})
$$
\end{small}where is ${\bf F}_s$ the symmetric part of ${\bf F}$, all
transformations $\bf{\Theta}$ corresponding to the same ${\bf M}$ lead
to the same eigenvalues for ${\bf F}_s \ $, and therefore to the same
contraction rate $| \lambda_{max} |$.

Consider now a hybrid case \cite{chemical}, consisting of a continuous
system $$ \dot{\bf x} = {\bf f}({\bf x}, t) $$ which is switched to a
discrete system $$ {\bf x}_{i+1} = {\bf f}_i({\bf x}_i, i) $$ every
$\Delta t_i$ for one discrete step. Letting, in the {\it same
coordinate system ${\bf \Theta}$}, $\bar{\lambda}$ be the largest
eigenvalue of the symmetric matrix ${\bf F}^T+{\bf F}$, and
$\bar{\lambda}_{i}$ be the largest eigenvalue of ${\bf F}_i^T {\bf
F}_i$ (the corresponding discrete-time quantity, where ${\bf F}_i
={\bf \Theta}_{i+1} \frac{\partial {\bf f}_i}{\partial {\bf x}_i} {\bf
\Theta}_i^{-1}$, see~\cite{chemical}), a sufficient condition for the
overall system to be contracting is
\begin{equation} \label{eq:hybrid_condition}
\exists \ \alpha < 1, \forall i,\ \ \ \ 0 \le
\bar{\lambda}_i e^{\bar{\lambda} \Delta t_i} \le \alpha
\end{equation}

Contraction theory proofs and this paper make extensive use of {\it
virtual displacements}, which are differential displacements at fixed
time borrowed from mathematical physics and optimization
theory. Formally, if we view the position of the system at time $t$ as
a smooth function of the initial condition ${\bf x}_o$ and of time, $\
{\bf x} = {\bf x}({\bf x}_o ,t)\ $, then $\ \delta {\bf x} =
\frac{\partial {\bf x}}{\partial {\bf x}_o} \ d {{\bf x}_o}\ $.

\end{document}